\documentclass[12pt]{amsart}
\usepackage{amsmath, amsthm, euscript}
\swapnumbers
\theoremstyle{plain}
\newtheorem{thm}{Theorem}[section]
\newtheorem{lem}[thm]{Lemma}
\newtheorem{prop}[thm]{Proposition}
\newtheorem{cor}[thm]{Corollary}

\theoremstyle{definition}
\newtheorem{setup}[thm]{Setup}
\newtheorem{exam}[thm]{Example}
\newtheorem{Weyl}[thm]{Example: Inverse Weyl Power Series}
\newtheorem*{Ack}{Acknowledgement}

\theoremstyle{remark}
\newtheorem{note}[thm]{}

\def\ann{\operatorname{ann}}
\def\gr{\operatorname{gr}}
\def\rgldim{\operatorname{rgldim}}
\def\rKdim{\operatorname{rKdim}}
\def\clKdim{\operatorname{clKdim}}
\def\Kdim{\operatorname{Kdim}}
\def\yi{y^{-1}}
\def\Yi{Y^{-1}}
\def\zi{z^{-1}}
\def\xi{x^{-1}}
\def\Xi{X^{-1}}
\def\x{\overline{x}}

\def\m{\mathfrak{m}}
\def\n{\mathfrak{n}}
\def\a{\mathfrak{a}}

\begin{document}

\title{Noetherian Skew Inverse Power Series Rings}

\author{Edward S. Letzter}

\author{Linhong Wang}

\address{Department of Mathematics\\
        Temple University\\
        Philadelphia, PA 19122-6094}
      
      \email[\emph{first author}]{letzter@temple.edu }

      \email[\emph{second author}]{lhwang@temple.edu}
      
      \thanks{Research of the first author supported in part by grants
        from the National Security Agency. Research of the first author
        also supported in part by Leverhulme Research Interchange Grant
        F/00158/X (UK).}



\begin{abstract} We study skew inverse power series extensions
  $R[[\yi;\tau,\delta]]$, where $R$ is a noetherian ring equipped with an
  automorphism $\tau$ and a $\tau$-derivation $\delta$. We find that these
  extensions share many of the well known features of commutative power series
  rings. As an application of our analysis, we see that the iterated skew
  inverse power series rings corresponding to $n$th Weyl algebras are
  complete, local, noetherian, Auslander regular domains whose right Krull
  dimension, global dimension, and classical Krull dimension are all equal to
  $2n$.
\end{abstract}

\maketitle


\section{Introduction} 

Let $R$ be a ring equipped with an automorphism $\tau$ and a
$\tau$-derivation $\delta$. The skew Laurent series ring
$R((y;\tau))$, when $\delta = 0$, and the pseudodifferential operator
ring $R[[\yi;\delta]]$, when $\tau$ is the identity, are well known,
classical objects. (See, e.g., \cite{Tug} for relevant history and new
results; pseudodifferential operator rings appear in \cite{Bjo-2} as
rings of germs of micro-local differential operators.)  These rings
provide noncommutative generalizations of commutative power series and
Laurent series rings. Other generalizations include the suitably
conditioned skew power series rings $R[[y;\tau,\delta]]$ recently
studied in \cite{SchVen},\cite{Ven}. In the present paper we study the
inverse skew power series rings $R[[\yi;\tau,\delta]]$, which turn out
to be particularly well behaved analogues of commutative power series
rings.  As an application, we provide (apparently) new examples of
complete, local, noetherian, Auslander regular domains corresponding
to $n$th Weyl algebras; working over a field, these iterated skew
inverse power series rings have right Krull dimension, global
dimension, and classical Krull dimension all equal to $2n$.

Our approach is largely derived from the commutative case and from the studies
of skew polynomial rings found in \cite{Goo} and \cite{GooLet}. Also playing
key roles are well known filtered-graded arguments (cf., e.g.,
\cite{LiVOy},\cite{NasVOy}).

The paper is organized as follows: Section 2 examines the interplay between
the ideal structure of $R$ and $R[[\yi;\tau,\delta]]$, focusing on primality,
locality, completeness, dimensions, and Auslander regularity. Section 3
considers iterated extensions, including the examples derived from Weyl
algebras.

\begin{Ack} We are grateful to K. A. Brown for helpful remarks. We are
  also grateful to the referee for careful reading, detailed comments,
  and helpful suggestions.
\end{Ack}

\section{Skew Inverse Power Series Rings} 

In this section we present some basic properties of skew inverse power series
extensions. 

\begin{setup} \label{2setup} We first set notation (to remain in effect
  for the remainder of this paper) and briefly review the constructions basic
  to our study; the reader is referred to \cite{Coh} and \cite{GooWar} for
  details.
  
  (i) To start, $R$ will denote an associative unital ring, equipped with a
  ring automorphism $\tau$ and a (left) $\tau$-derivation $\delta$. In other
  words, $\delta{\colon}R \rightarrow R$ is an additive map for which
  $\delta(ab) = \tau(a)\delta(b) + \delta(a)b$ for all $a, b \in R$. An
  element $\lambda$ of $R$ is \emph{$\tau$-$\delta$-scalar} if $\tau(\lambda)
  = \lambda$ and $\delta(\lambda) = 0$. We also refer to the \emph{(left) skew
    derivation} $(\tau,\delta)$ on $R$.

  Recall that the skew polynomial ring $R[y;\tau, \delta]$, comprised of
  polynomials
\[ r_n y^n + r_{n-1} y^{n-1} + \cdots + r_o, \]
for $n = 0,1,2,\ldots$ and $r_0,\ldots,r_n \in R$, is constructed via the
multiplication rule
\[yr = \tau(r)y + \delta(r), \]
for all $r \in R$. 

(ii) Set $S = R[[\yi;\tau,\delta]]$, the ring of formal skew power series in
$\yi$,
\[ \sum_{i=0}^\infty y^{-i} r_i , \]
for $r_0,r_1,\ldots \in R$. Multiplication is determined by the rule
\[ r\yi = \sum _{i=1}^\infty y^{-i}\tau\delta^{i-1}(r), \]
for $r \in R$, derived from
\[ r\yi = \yi \tau(r) + \yi \delta(r) \yi .\]

(iii) Since $\tau$ is an automorphism, we have
\[ \yi r = \tau^{-1}(r)\yi - \yi\delta \tau^{-1}(r)\yi ,\]
from which we can deduce the formula
\[ \yi r = \sum_{i=1}^\infty \tau^{-1}\big( -\delta \tau^{-1} \big)^{i-1}(r) y^{-i} .\]
Consequently, we can write coefficients of power series in $S$ on either the
right or left.

(iv) It follows from (ii) and (iii) that $\yi$ is normal in $S$ (i.e., $S\yi =
\yi S$). Localizing $S$ at the powers of $\yi$, we obtain $S' :=
R((\yi;\tau,\delta))$, the ring of formal skew Laurent series in $\yi$,
\[ \sum_{i=-n}^\infty y^{-i} r_i , \]
for $r_0,r_1,\ldots \in R$ and non-negative integers $n$. In view of (iii) we
may write coefficients of Laurent series in $S'$ on either the left or right.

Viewing the preceding power and Laurent series, with coefficients on the
right, as $\infty$-tuples, we see that $S$ and $S'$ are naturally isomorphic
as right $R$-modules to infinite direct products of copies of $R$. We
similarly obtain left direct product structures for $S$ and $S'$.

(v) To remind us that $\yi$ is not a unit in $S$, we will often use the
substitution $z := \yi$. Note that $z$ is a regular element of $S$ and that
there is a natural isomorphism of $R$ onto $S/\langle z \rangle$.

(vi) Let $f = r_0 + zr_1 + z^2r_2 + \cdots$ be a nonzero power series in $S$,
for $r_0,r_1,\ldots \in R$. The \emph{right initial coefficient} for $f$ will
be the first appearing nonzero $r_i$. Writing coefficients on the left we can
similarly define the \emph{left initial coefficient} of $f$. We will refer to
$r_0$ as the \emph{constant coefficient} of $f$. Note that the constant
coefficient is the same whether we use left or right coefficients. We identify
$R$ with the subring of $S$ consisting of power series all of whose
nonconstant coefficients are zero.

If $R$ is a domain then an easy argument, involving initial coefficients of
products of power series in $S$, shows that $S$ must also be a domain.

\end{setup}

We now consider ideals, particularly prime ideals and generalizations.

\begin{note} \label{2r1} Let $J$ be an ideal of $S$. The
  set comprised of $0$ together with the right initial coefficients of power
  series in $J$ forms an ideal of $R$, as does the set comprised of $0$
  together with the left initial coefficients. Also, the set of constant
  coefficients of $J$ forms an ideal of $R$.
\end{note}

\begin{note} We will also make frequent use of the following terminology: An
  ideal $I$ of $R$ for which $\tau(I) \subseteq I$ is a \emph{$\tau$-ideal},
  and a $\tau$-ideal $I$ of $R$ for which $\delta(I) \subseteq I$ is a
  \emph{$\tau$-$\delta$-ideal}. A $\tau$-$\delta$-ideal $P$ of $R$ such that
  $AB \subseteq P$ only if $A \subseteq P$ or $B \subseteq P$, for all
  $\tau$-$\delta$-ideals $A$ and $B$ of $R$, is \emph{$\tau$-$\delta$-prime}.
  We will say that $R$ is $\tau$-$\delta$-prime if $0$ is a
  $\tau$-$\delta$-prime ideal. We similarly define \emph{$\tau$-prime} ideals
  and say that $R$ is $\tau$-prime when $0$ is a $\tau$-prime ideal. The
  reader is referred to \cite{GooLet} for background, particularly on the
  interplay between the ideal theory in $R[y;\tau,\delta]$ and the
  $\tau$-$\delta$-ideal theory of $R$. 

  Note, when $R$ is right or left noetherian, that an ideal $I$ of $R$
  is a $\tau$-ideal if and only if $\tau(I) = I$.
\end{note}

\begin{note} \label{2r2} Suppose that $R$ is right or left noetherian,
and let $I$ be an ideal of $R$. It follows from \cite[Remarks $4^*$, $5^*$,
p$.$ 338]{GolMic} that $I$ is $\tau$-prime if and only if $I = Q \cap \tau(Q)
\cap \cdots \cap \tau^t (Q)$ for some prime ideal $Q$ of $R$ and some
non-negative integer $t$ such that $\tau^{t+1}(Q) = Q$.
\end{note}

\begin{note} \label{2r3} (i) Consider the filtration
\[ S = \langle z \rangle ^0 \supset \langle z \rangle^1 \supset \langle z
\rangle^2 \supset \cdots , \]
and the induced $\langle z \rangle$-adic topology on $S$. The filtration is
exhaustive (i.e., $S = \langle z \rangle ^0$), separated (i.e., the $\langle z
\rangle$-adic topology is Hausdorff), and complete (i.e., Cauchy sequences
converge in the $\langle z \rangle$-adic topology); see, e.g., \cite[Chapter
D]{NasVOy} for background.  Moreover, the associated graded ring,
\[ \gr(S) \; = \; R \oplus \langle z \rangle/ \langle z \rangle ^2 \oplus
\langle z \rangle ^2/\langle z \rangle^3 \oplus \cdots, \]
is isomorphic to the skew polynomial ring $R[x;\tau^{-1}]$. 

(ii) It follows from (i) and (\ref{2setup}v) that $S$ is left or right
noetherian if and only if the same holds for $R$ (see, e.g.,
\cite[D.IV.5]{NasVOy}).

(iii) Suppose that $R$ is right noetherian. It follows from (i) and \cite[p$.$
87, Proposition]{LiVOy} that the $\langle z \rangle$-adic filtration on $S$
is \emph{Zariskian} in the sense of \cite{LiVOy} (cf$.$ \cite{Bjo-1}).
\end{note}

The proof of the next lemma is omitted.

\begin{lem} \label{2t1} Assume that $R$ is right or left noetherian, and let
  $I$ be a $\tau$-$\delta$-ideal of $R$.  {\rm (i)} $IS = SI = SIS$.

{\rm (ii)}
\[ \left. \left\{ \sum_{i=0}^\infty z^ia_i \; \right| \; a_0,a_1,\ldots
  \in I \right\} \quad = \quad \left. \left\{ \sum_{i=0}^\infty b_iz^i \;
  \right| \; b_0,b_1,\ldots \in I \right\}\]
is a two-sided ideal of $S$.
\end{lem}

\begin{note}
 The ideal of $S$ described in (ii) of the preceding lemma will
be denoted by $I[[\yi;\tau,\delta]]$ or $I\langle\langle z \rangle\rangle$.
Note that
\[S/I\langle\langle z\rangle\rangle \; = \;
R[[\yi;\tau,\delta]]\big/I[[\yi;\tau,\delta]] \; \cong \;
\left(R/I\right)[[\yi;\tau,\delta]] \; \cong \; \left(R/I\right)\langle\langle
z \rangle\rangle.\]
Also note that $I\langle\langle z \rangle\rangle \cap R = I$.
\end{note}

The following and its proof are adapted from standard commutative arguments.

\begin{prop} \label{2t2} Assume that $R$ is right or left noetherian,
  and let $I$ be a $\tau$-$\delta$-ideal of $R$. Then $I\langle\langle z \rangle\rangle = IS =
  SI$.
\end{prop}

\begin{proof} First suppose that $R$ is left noetherian, and so $I =
  Ra_1 + \cdots + Ra_n$, for $a_1,\ldots,a_n \in I$. Choose
\[ \sum_{i=0}^\infty z^i b_i \; \in \; I\langle\langle z \rangle\rangle ,\]
with $b_0,b_1,\ldots \in I$. Then, for suitable choices of $r_{ij} \in R$,
\[ \sum_{i=0}^\infty z^ib_i \; = \; \sum_{i=0}^\infty z^i(r_{1i}a_1 + \cdots +
r_{ni}a_n) \; = \; \left(\sum_{i=0}^\infty z^ir_{1i}\right)a_1 + \cdots +
\left(\sum_{i=0}^\infty z^ir_{ni}\right)a_n \; \in \; SI .\]
Hence $I\langle\langle z \rangle\rangle \subseteq IS = SI$. Of course, $IS
\subseteq I\langle\langle z \rangle\rangle$, and the proposition follows in
this case. The case when $R$ is right noetherian follows
similarly, writing coefficients of power series on the left.
\end{proof}

\begin{note} \label{extends} An important special case occurs when
  $\tau$ can be extended to automorphisms of $S$ and $S'$. This
  situation occurs, for instance, when $\delta\tau = \tau\delta$ (as
  operators on $R$); in this case we can extend $\tau$ to $S$ and $S'$
  by setting $\tau(y^{\pm 1}) = y^{\pm 1}$. (A proof of this assertion
  will follow from the next paragraph, setting $q = 1$.) Note that
  $\delta$ and $\tau$ satisfy the equation $\delta\tau = \tau\delta$
  when $\tau$ is the identity or when $\delta = 0$.
  
  More generally, suppose (temporarily) that $\delta\tau = q\tau\delta$ for
  some central unit $q$ of $R$ such that $\tau(q) = q$ and $\delta(q) = 0$.
  (See, e.g., \cite{Goo},\cite{GooLet}.) Observe, in view of (\ref{2setup}ii),
  that
\[ \tau(r)q\yi = \sum _{i=1}^\infty qy^{-i}\tau\delta^{i-1}\tau(r) =
\sum_{i=1}^\infty (q^i y^{-i})\tau(\tau\delta^{i-1}(r)) .\]
It follows that $\tau$ is compatible with multiplication in $S$ and so
extends to an automorphism of $S$, with $\tau(\yi) = q\yi$.  It also
follows that $\tau$ extends from $S$ to an automorphism of $S'$, with
$\tau(y) = q^{-1}y$.

Even more generally, removing the assumption that $\delta\tau =
q\tau\delta$, we see whenever $\tau$ extends to an automorphism of $S$
for which $\tau(\yi) = \mu\yi$, for a unit $\mu$ of $S$, that $\tau$
extends from $S$ to an automorphism of $S'$.

\end{note}

\begin{prop} \label{2t3} {\rm (i)} Assume that $R$ is right or left
  noetherian. Let $P$ be a prime ideal of $S$, and suppose that $P\cap
  R$ is a $\tau$-$\delta$-ideal of $R$. Then $P\cap R$ is
  $\tau$-$\delta$-prime. In particular, if $S$ is prime then $R$ is
  $\tau$-$\delta$-prime. {\rm (ii)} Suppose that $\tau$ extends to an
  automorphism of $S$, that $P$ is a $\tau$-prime ideal of $S$, and
  that $P\cap R$ is a $\tau$-$\delta$-ideal of $R$. Then $P\cap R$ is
  a $\tau$-$\delta$-prime ideal of $R$. In particular, if $S$ is
  $\tau$-prime then $R$ is $\tau$-$\delta$-prime.
\end{prop}

\begin{proof} We prove (ii); the proof of (i) is similar and easier. Let $I$ and $J$ be $\tau$-$\delta$-ideals of $R$ such that $IJ \subseteq P\cap R$.
  By (\ref{2t1}i),
\[ (SIS)(SJS) = SIJ \subseteq P. \]
Note that $SIS$ and $SJS$ are $\tau$-ideals of $S$, and so either 
$SIS \subseteq P$ or $SJS \subseteq P$. Therefore, either $I \subseteq
(SIS)\cap R \subseteq P\cap R$ or $J \subseteq (SJS)\cap R \subseteq P\cap R$.
\end{proof}

\begin{note} Let $I$ be a $\tau$-$\delta$-ideal of $R$. It follows from (i),
  in the preceding proposition, that if $I\langle\langle z \rangle\rangle$ is
  prime then $I$ is $\tau$-$\delta$-prime. When $\tau$ extends to an
  automorphism of $S$, it follows from (ii) that if $I\langle\langle z
  \rangle\rangle$ is $\tau$-prime then $I$ is $\tau$-$\delta$-prime.
\end{note}

The next example shows that $S$ need not be
prime -- or, indeed, semiprime -- when $R$ is $\tau$-$\delta$-prime.

\begin{exam} \label{2ex} We follow \cite[2.8]{Goo} and
  \cite[3.1]{GooLet}, where detailed calculations justifying our assertions
  can be found. To start, let $k$ be a field, and let $\alpha$ be the
  automorphism of the ring $k^4$ given by $\alpha(a,b,c,d) = (b,c,d,a)$. Set
  $U = k^4[x;\alpha]$ and $T = k^4[x^{\pm 1}; \alpha]$. Extend $\alpha$ to an
  automorphism of $T$ by setting $\alpha(x) = x$, and let $\tau =
  \alpha^{-1}$.  Let $\delta$ denote the following $\tau$-derivation of
  $T$:
\[ \delta(t) = (0,0,0,1)\xi t - \tau(t)(0,0,0,1)\xi,\]
for $t \in T$. Then $\delta$ restricts to a $\tau$-derivation of $U$, and
$Ux^4$ is a maximal proper $\tau$-$\delta$-ideal of $U$. Now set $R = U/Ux^4 =
k^4\langle \x \rangle$, where $\x$ denotes the image of $x$ in $R$, and set $v
= (1,0,0,0)\x \in R$. Then $R$ contains no $\tau$-$\delta$-ideals other than
$0$ and $R$, and so in particular $R$ is $\tau$-$\delta$-prime.

As in \cite[p$.$ 16]{GooLet}, $v{\x}^2 \ne 0$. Also,
since $\x$ is normal in $R$, $v{\x}^2$ is normal in $R$. Again as in
\cite[p$.$ 16]{GooLet},
\[ y^4 v = v y^4 \]
in $R[y;\tau,\delta]$, and
\[ v{\x}^2y^tv = 0 ,\]
for $t = 0,1,2,\ldots$. Now set $S = R[[y^{-1};\tau,\delta]].$ In $S$, by
above,
\[ y^{-4}v = v y^{-4} .\]
Therefore, for all $r \in R$, all $l = 0,1,2,\ldots$, and all $t =
0,1,2,\ldots$,
\[ v{\x}^2ry^{-4l+t}v{\x}^2 \; = \; r'(v{\x}^2y^tv)y^{-4l}{\x}^2 \; = \; 0 
\]
for some $r' \in R$. Hence, $v{\x}^2Sv{\x}^2 = 0$, and $S$ is not semiprime.
\end{exam}

\begin{note} \label{2r4} We now return to $S' = R((y^{-1};\tau,\delta))$, 
the localization of $S$ at powers of $z$, recalled in (\ref{2setup}iv).
Suppose (for now) that $R$ is right or left noetherian.

(i) There is a natural bijection
\[ \left\{ \text{semiprime ideals of $S$ not containing $z$} \right\}
\quad \longleftrightarrow \quad \left\{ \text{semiprime ideals of $S'$}
\right\}, \]
obtained via the extension map $I \mapsto IS'$, for semiprime ideals
$I$ of $S$ not containing $z$, and the contraction map $J \mapsto
J\cap S$, for semiprime ideals $J$ of $S'$. Also, if $I$ is a prime
ideal of $S$ not containing $z$ then $IS'$ is a prime ideal of $S'$,
and if $J$ is a prime ideal of $S'$ then $J \cap S$ is a prime ideal
of $S$. (See, e.g., \cite[10.18]{GooWar}.)

(ii) Suppose, for the moment, that $\tau$ extends to an automorphism
of $S$ and further extends from $S$ to an automorphism of $S'$; this
situation was considered in (\ref{extends}). Using (i) and
(\ref{2r2}), we see in this case that extension and contraction
provide a natural bijection between the set of $\tau$-prime ideals of
$S$ not containing $z$ and the set of $\tau$-prime ideals of $S'$.
\end{note}

The following and its proof are directly adapted from \cite[3.2]{Goo} and
\cite[3.3]{GooLet}. We see that the behavior exhibited in (\ref{2ex}) cannot
occur when $\tau$ extends to a automorphisms of $S$ and $S'$.

\begin{prop} \label{2t4} Assume that $R$ is $\tau$-$\delta$-prime and
  that $R$ is right or left noetherian. Also assume that $\tau$
  extends to an automorphism of $S$ and further extends from $S$ to an
  automorphism of $S'$. {\rm (i)} $S'$ is $\tau$-prime.  {\rm (ii)}
  $S$ is $\tau$-prime.
\end{prop}

\begin{proof} (i) Assume that $S'$ is not $\tau$-prime. Then there
  exist $\tau$-ideals $I$ and $J$ of $S'$, both nonzero, such that $IJ
  = 0$. Without loss of generality, $I$ is the left annihilator in
  $S'$ of $J$, and $J$ is the right annihilator in $S'$ of $I$. (This
  simplification makes use of the fact that $\tau(I) = I$ and $\tau(J)
  = J$.)
  
  Let
\[ f \; = \; \sum_{i=j}^\infty z^i a_i \]
be a nonzero power series in $I$, with right initial coefficient $a = a_j$,
and let
\[ g \; = \; \sum_{i=k}^\infty b_i z^i \]
be a nonzero power series in $J$, with coefficients written on the left.  Of
course, $fg = 0$. Now suppose that $ag \ne 0$. Then we can choose $l$ minimal
such that $ab_l \ne 0$. Therefore, $\tau^l(ab_l)$ is the (nonzero) leading
right coefficient of $fg$, contradicting the fact that $fg = 0$. Hence, $ag =
0$, and
\[ 0 \ne a \in \left(\ann_{S'}J \right)\cap R = I \cap R .\]
In particular, $I \cap R \ne 0$. Similar reasoning shows that $J\cap R \ne 0$.

Now choose $r \in I \cap R$.  Observe that 
\[ \delta(r) = yr - \tau(r)y = \zi r - \tau(r)\zi \in I, \]
and so $\delta(r) \in I \cap R$. Thus $I\cap R$ is a nonzero
$\tau$-$\delta$-ideal of $R$. Similarly, $J\cap R$ is a nonzero 
$\tau$-$\delta$-ideal of $R$. But $(I\cap R)(J\cap R) = 0$, contradicting our
assumption that $R$ is $\tau$-$\delta$-prime. Part (i) follows.

(ii) This follows from (i) and (\ref{2r4}ii).
\end{proof}

Next, we consider maximal ideals, the Jacobson radical, and locality.

\begin{lem} \label{2t5} {\rm (i)} Let $f = 1 + zr_1 + z^2r_2 + \cdots
  \in S$ for $r_1,r_2,\ldots \in R$. Then $f$ is a unit in $S$. {\rm
    (ii)} Let $g \in S$. If the constant coefficient of $g$ is a unit
  in $R$, then $g$ is a unit in $S$. {\rm (iii)} $z$ is contained in
  the Jacobson radical $J(S)$ of $S$.  
\end{lem}

\begin{proof} Set $u = 1 + za_1 + z^2a_2 + \cdots \in S$, for
  $a_1,a_2,\ldots \in R$, and set $v = 1 + zb_1 + z^2b_2 + \cdots$, for
  $b_1,b_2,\ldots \in R$. Then
\[ uv \; = \; 1 + z(a_1 + b_1) + z^2(a_2 + b_2 + p_2(a_1,b_1)) + z^3(a_3 +
  b_3 + p_3(a_1,a_2,b_1,b_2)) + \cdots ,\]
  where $p_i(a_1,\ldots,a_{i-1},b_1,\ldots,b_{i-1}) \in R$ depends only on
  $a_1,\ldots,a_{i-1}$ and $b_1,\ldots,b_{i-1}$, for $i = 2,3,\ldots$.

  If $a_1,a_2,\ldots$ are arbitrary then we can choose
  $b_1,b_2,\ldots$ such that $uv = 1$, and if $b_1,b_2,\ldots$ are
  arbitrary then we can choose $a_1,a_2,\ldots$ such that $uv =
  1$. Part (i) now follows, and part (ii) follows easily from (i). To
  prove part (iii), let $a$ be an arbitrary element of $S$. By (i), $1
  + za$ is a unit in $S$, and so $z$ is contained in $J(S)$. 
\end{proof}

\begin{prop} \label{2t6} {\rm (i)} $J(S) = J(R) + \langle z \rangle$. {\rm
    (ii)} Let $P$ be a maximal ideal of $S$. Then $P = Q + \langle z \rangle$
  for some maximal ideal $Q$ of $R$. In particular, if $R$ has a unique
  maximal ideal $Q$, then $Q + \langle z \rangle$ is the unique maximal ideal
  of $S$. {\rm (iii)} Suppose that $R$ has a unique maximal ideal $Q$, and
  suppose further that every element of $R$ not contained in $Q$ is a unit.
  Then every element of $S$ not contained in $Q + \langle z \rangle$ is a
  unit.
\end{prop}

\begin{proof} (i) By (\ref{2t5}iii), $z$ is contained in $J(S)$. It then
  follows from the natural isomorphism of $R$ onto $S/\langle z \rangle$ that
  $J(S) = J(R) + \langle z \rangle$.

  (ii) Let $I$ be the ideal in $R$ of constant coefficients of power series in
  $P$. Suppose $I$ is not contained in a maximal ideal of $R$.  Then $1 + zr_1
  + z^2r_2 + \cdots \in P$, for some choice of $r_1,r_2,\ldots \in R$, and so
  $P = S$ by (\ref{2t5}i). Therefore, $I$ is contained in some maximal ideal
  $Q$ of $R$, and so $P \subseteq Q + \langle z \rangle$.  Since $P$ is a
  maximal ideal of $S$, and since $Q + \langle z \rangle$ is also a maximal
  ideal of $S$, we see that $P = Q + \langle z \rangle$.

  (iii) If $f \notin Q + \langle z \rangle$ then the constant coefficient of
  $f$ is a unit, and so $f$ is a unit in $S$ by (\ref{2t5}ii).
\end{proof}

Now we turn to topological properties.

\begin{thm} \label{2t7} Assume that $R$ is right or left
  noetherian. Let $I$ be a $\tau$-$\delta$-ideal of $R$, and suppose
  that the $I$-adic filtration of $R$ is separated and complete. Set
  $J = I + \langle z \rangle$. Then the $J$-adic filtration of $S$ is
  separated and complete.
\end{thm}

\begin{proof} First note, for positive integers $t$, that
\[J^t \; = \; z^t S + z^{t-1}SI + \cdots + zSI^{t -1} + SI^t ,\]
by (\ref{2t1}i). Hence, given $r_0 + zr_1 + z^2r_2 + \cdots \in J^t$,
for $r_0,r_1,\ldots \in R$, it follows that:
\[\text{$r_j \in I^{t-j}$ for all $j = 0,1,2,\ldots$} \quad \tag{$\ast$}\]

Now, to prove separation, let
\[ b \; = \; a_0 + za_1 + z^2a_2 + \cdots \quad \in \quad \bigcap_{i=1}^\infty J^i ,\]
for $a_0,a_1,\ldots \in R$, and suppose that $b \ne 0$. Let $j$ be a
non-negative integer minimal such that $a_j \ne 0$. Then, for any $t >
j$, it follows from ($\ast$) that
\[a_j \; \in \; J^t \; = \; z^t S + z^{t-1}SI + \cdots + z^j SI^{t-j} + \cdots
+ zSI^{t-1} + SI^t .\]

It follows that
\[ a_j \; \in \; \bigcap_{l=1}^\infty I^l \; = \; 0, \]
a contradiction. Therefore, $b = 0$, and the $J$-adic filtration on $S$ is
separated.

Next, let
\[\begin{aligned} s_1 \; &= \; r_{10} + zr_{11} + z^2r_{12} + \cdots, \\
  s_2 \; &= \; r_{20} + zr_{21} + z^2r_{22} + \cdots, \cdots
\end{aligned}\]
be a Cauchy sequence in $S$ with respect to the $J$-adic topology. For each
non-negative integer $m$ we obtain a sequence $\{ r_{nm}\}_{n=1}^\infty$, and
for each integer $n \geq m$ it follows that $s_u - s_v \in J^n$ for
sufficiently large $u$ and $v$. Note that the $m$th coefficient of $s_u - s_v$
is $r_{um} - r_{vm}$. It follows from ($\ast$) that $r_{um}
- r_{vm} \in I^{n-m}$. Therefore, $\{ r_{nm}\}_{n=1}^\infty$ is a Cauchy
sequence in $R$ with respect to the $I$-adic topology. Set
\[ r_m = \lim_{n\rightarrow \infty} r_{nm}, \quad \text{and} \quad s = r_0 +
zr_1 + z^2r_2 + \cdots \; .\]
Again using $(\ast)$, it follows that the sequence $s_1, s_2, \ldots$
converges to $s$ in the $J$-adic topology, and so $S$ is complete.
\end{proof}

Recall that a ring is \emph{local\/} provided its Jacobson radical is a
coartinian maximal ideal. By a \emph{complete local} ring we will always mean
a local ring, with Jacobson radical $J$, such that the $J$-adic filtration
is separated and complete

\begin{cor} \label{2t8} Suppose that $R$ is a complete local ring with
  unique primitive ideal $\m$, and suppose that $\delta(\m) \subseteq
  \m$. Then $S$ is a complete local ring whose unique primitive ideal
  is $\m + \langle z \rangle$.
\end{cor}

\begin{proof} Set $J = \m + \langle z \rangle$. It follows from
  (\ref{2t6}i) that $J$ is the Jacobson radical of $S$, and it follows
  from (\ref{2t6}ii) that $J$ is the unique maximal (and also unique
  primitive) ideal of $S$. Next, $\m$ must be $\tau$-stable, and so
  $\m$ is a $\tau$-$\delta$-ideal. That $S$ is a complete local ring
  now follows from (\ref{2t7}).
\end{proof}

We now consider dimensions and Auslander regularity.

\begin{note} \label{2r5} Assume that $R$ is right noetherian, and let $n$ be a
  non-negative integer. Recall from (\ref{2setup}iv, v) that $z$ is a normal
  regular element of $S$ and from (\ref{2t5}iii) that $z$ is contained in the
  Jacobson radical of $S$.

  (i) Suppose that that $\rKdim R$, the right Krull dimension of $R$ (see,
  e.g., \cite{GooWar}), is equal to $n$. It follows from \cite[1.8]{Wal} that
  $\rKdim S = n+1$.
  
  (ii) Suppose that $\clKdim R$, the classical Krull dimension of $R$ (see,
  e.g., \cite{GooWar}), is equal to $n$. Since $z$ is regular in $S$, it
  follows (e.g., from \cite[11.8]{GooWar}) that $z$ is regular modulo the
  prime radical of $S$.  Therefore, no prime ideal of $S$ containing $z$ can
  be a minimal prime ideal of $S$. Consequently, the classical Krull dimension
  of $S$ is strictly greater than the classical Krull dimension of $S/zS \cong
  R$; in other words, $\clKdim S \geq n + 1$.

  Now recall that $\clKdim S \leq \rKdim S$. In particular, by (i), if
  $\clKdim R = \rKdim R = n$ then $\clKdim S = \rKdim S =
  n+1$. Furthermore, when $J(S)$ has a normalizing set of generators
  (a condition satisfied by the examples considered in \S
  \ref{iterated}), it can also be deduced from \cite[2.7]{Wal} that
  $\clKdim S = \rKdim S$.
  
  (iii) When the right global dimension $\rgldim R$ is equal to $n$ it follows
  from \cite[1.3]{Wal} that $\rgldim S = n+1$.

  (iv) Assume now that $R$ is noetherian (on both sides). Recall from
  (\ref{2r3}i) that the $\langle z \rangle$-adic filtration on $S$ is
  exhaustive, separated, and complete, with $\gr(S) \cong R[x;\tau^{-1}]$.
  Also recall, from (\ref{2r3}iii), that this filtration is Zariskian. Now
  suppose further that $R$ is Auslander regular (see, e.g., \cite[\S
  III.2]{LiVOy}). It then follows from \cite[p$.$ 174, Theorem]{LiVOy} that
  $\gr(S)$ is Auslander regular. Therefore, by \cite[p$.$ 152,
  Theorem]{LiVOy}, $S$ is Auslander regular. (The preceding can also be
  deduced from \cite{Bjo-1}.)

\end{note}

\section{Iterated Skew Inverse Power Series Rings} \label{iterated}

We now apply the analysis of the preceding section to iterative constructions.

\begin{setup} \label{3setup} Let $C$ be a commutative, complete, regular,
  local, noetherian domain with maximal ideal $\m$ and residue field $k$. Let
  $m$ be a positive integer, and set $R_0 = C$.  For each $1 \leq i \leq m$
  let
\[ R_i = C [[\yi_1 ; \tau_1, \delta_1]]
\cdots[[\yi_i; \tau_i, \delta_i]] ,\]
where $\tau_i$ is a $C$-algebra automorphism of $R_{i-1}$, and where
$\delta_i$ is a left $C$-linear $\tau_i$-derivation of $R_{i-1}$. (Note,
then, that $\tau_1$ is the identity and $\delta_1$ is trivial.) Set $A = R_m$,
and set $z_i = \yi_i$ for each $i$. Let $\n_i$ denote the ideal of
$R_i$ generated by $\m$ and $z_1,\ldots,z_i$. Set $\n = \n_m$ and $\n_0 = \m$.
\end{setup}

\begin{note} \label{3r2} Fix $1 \leq i \leq m$. It follows from
  (\ref{2setup}iv), for $1 \leq j \leq i$, that
\[ R_iz_i + R_iz_{i-1} + \cdots + R_iz_j \; = \; z_iR_i + z_{i-1}R_i +
\cdots + z_jR_i \]
and, for $1 < j \leq i$, that $z_{j-1}$ is normal in $R_i$ modulo the ideal
$R_iz_i + \cdots + R_iz_j$.  Furthermore, since $C \cong
R_i/\langle z_1,\ldots,z_i\rangle$, we see that
\[ \n_i \; = \; R_i\m + R_iz_1 + R_iz_2 + \cdots + R_iz_i \; = \; 
\m R_i + z_1R_i + z_2R_i + \cdots + z_iR_i . \]
In particular, since $\m$ is a finitely generated ideal of $C$, it follows
that $\n_i$ has a normalizing sequence of generators in $R_i$.

Observe, further, that $R_i/\n_i \cong k$.
\end{note}

We can now apply results and observations from the preceding section, as
follows.

\begin{prop} \label{3t1} Let $A = C[[\yi_1 ; \tau_1, \delta_1]] \cdots[[\yi_m;
  \tau_m, \delta_m]]$, as above. {\rm (i)} $A$ is a (left and right)
  noetherian, Auslander regular domain. {\rm (ii)} $\rgldim A = \clKdim A =
  \rKdim A = \Kdim C + m$. {\rm (iii)} Every element of $A$ not contained in
  $\n$ is a unit, and $\n$ is the unique maximal ideal of $A$. {\rm (iv)}
  $J(A) = \n$. {\rm (v)} Suppose, for all $1 \leq i \leq m$, that
  $\delta_i(\n_{i-1}) \subseteq \n_{i-1}$. Then $A$ is a complete local ring
  with unique primitive ideal $\n$.
\end{prop}

\begin{proof} (i) Noetherianity follows (inductively) from (\ref{2setup}vi)
  and (\ref{2r3}ii). Since $C$ is Auslander regular, the Auslander regularity
  of $A$ follows from (\ref{2r5}iv).

  (ii) Follows from (\ref{2r5}i--iii).

(iii) Follows from repeated applications of (\ref{2t6}iii).

  (iv) Follows from repeated applications of (\ref{2t6}i).

(v) Follows from repeated applications of
  (\ref{2t8}).
\end{proof}

\begin{Weyl} We first review the (slightly) quantized Weyl algebras considered
  in \cite[2.8]{GooLet} (cf$.$ references cited therein). To start, let $T$ be
  a ring equipped with an automorphism $\sigma$, and let $q$ be a central
  $\sigma$-scalar unit in $T$. Extend $\sigma$ to the unique automorphism of
  $T[X;\sigma^{-1}]$ for which $\sigma(X) = qX$, and let $d$ be an arbitrary
  central element of $T$. Following \cite[2.8]{GooLet}, there is a unique
  (left) $\sigma$-derivation $D$ on $T[X;\sigma^{-1}]$ such that $D(T) = 0$
  and $D(X) = d$. We obtain the ``quantized Weyl algebra''
  $T[X;\sigma^{-1}][Y;\sigma,D]$ with coefficients in $T$. Note that $YX = qXY
  + d$ and that
\[ D(X^i) \; = \; d(q^{i-1}+q^{i-2}+\cdots+1)X^{i-1}, \]
for all positive integers $i$. When $\sigma$ is the identity and $q = d =1$, we
have the usual Weyl algebra with coefficients in $T$.
  
Using \cite[1.3]{Goo}, we can extend $(\sigma,D)$ uniquely to a skew
derivation of $T[X^{\pm 1};\sigma^{-1}]$, with
\[\sigma(tX^{-i}) \; = \; q^{-i}\sigma(t)X^{-i} \quad \text{and} \quad
D(tX^{-i}) \; = \; -d(q^{-1}+q^{-2}+\cdots+q^{-i})\sigma(t)X^{-i-1} ,\]
for all positive integers $i$ and all $t \in T$. Note that $(\sigma,D)$
restricts to a skew derivation of $T[\Xi;\sigma^{-1}]$. Next, for
\[ f = t_0 + t_1X^{-1} + t_2X^{-2} +
\cdots \; \in \; T[[\Xi;\sigma^{-1}]] ,\] 
with $t_0,t_1,t_2,\ldots \in T$, set
\[ \sigma(f) \; = \; \sum_{i=0}^\infty \sigma(t_iX^{-i}) \; = \;
\lim_{n\rightarrow \infty} \sigma\left( \sum_{i=0}^n t_iX^{-i}\right),\]
and
\[ D(f) \; = \; \sum_{i=0}^\infty D(t_iX^{-i}) \; = \; 
\lim_{n\rightarrow \infty} D\left( \sum_{i=0}^n t_iX^{-i}\right),\]
where limits are taken with respect to the $\langle \Xi \rangle$-adic topology
on $T[[\Xi;\sigma^{-1}]]$. It is now not hard to check that $\sigma$ and $D$
define a skew derivation of $T[[\Xi;\sigma^{-1}]]$, and we can construct the
skew inverse power series ring $T[[\Xi;\sigma^{-1}]][[\Yi;\sigma,D]]$.

We now iterate a simplified version of the preceding. Continue to let $C$ be a
commutative, complete, regular, local noetherian domain. Having constructed
\[ T \; = C[[\Xi_1]][[\Yi_1;\sigma_1,D_1]]\cdots[[\Xi_{n-1}]]
[[\Yi_{n-1};\sigma_{n-1},D_{n-1}]], \]
for positive integers $n$ (with $T = C$ when $n=1$), we can construct
\[ W \; = \; C[[\Xi_1]][[\Yi_1;\sigma_1,D_1]]\cdots[[\Xi_n]]
[[\Yi_n;\sigma_n,D_n]] ,\]
using the above procedure, with $\sigma_n = \sigma$ equal to the identity on
$T$, with $q_n = q$ equal to some unit in $C$, and with $d_n = d$ equal to
some element of $C$.

Note that $W$ satisfies all of the hypotheses of (\ref{3t1}), parts (1) and
(2). In particular, $W$ is a noetherian, complete, local, Auslander regular
domain.

When $C$ is a field, the following also follow from (\ref{3t1}): The ideal
\[\a = \langle \Xi_1,\Yi_1,\ldots,\Xi_n,\Yi_n \rangle \]
is the unique primitive ideal of $W$, every element of $W$ not contained in
$\a$ is a unit in $W$, and $\rgldim W = \rKdim W = \clKdim W = 2n$.
\end{Weyl}


\end{document}